\input amstex
\documentstyle{amsppt}

\topmatter

\title Weak$^*$ sequential closures in Banach space
theory and their applications
\endtitle

\rightheadtext{Weak$^*$ sequential closures}

\author M.I.Ostrovskii
\endauthor

\address  Department of Mathematics,
The Catholic University of America, 
Washington, D.C. 20064, USA
\endaddress

\email ostrovskii\@cua.edu\endemail



\endtopmatter

\document

\head  \S 1 Introduction
\endhead

Let $X$ be a (real or complex)
Banach space, its dual Banach space will be denoted by
$X^*$. 
We use standard notation and terminology of Banach space
theory, see J.Lindenstrauss and L.Tzafriri \cite{LT}.
By a {\it subspace} we mean a linear, but not
necessarily closed, subspace. We also assume some knowledge of general
topology and ordinal numbers, see P.S.Aleksandrov\cite{A}. 

\definition{Definition 1.1}  Let $A$ be a subset of $X^*$.
The set of all limits of weak$^*$-convergent sequences
in $A$ is called the {\it weak$^*$ sequential closure} of $A$ and is
denoted by $A_{(1)}$. 
\enddefinition

S.Banach asked the following question (see \cite{Maz}).
\medskip

{\bf Question.} Let $X$ be separable Banach space and $A$
be a subspace of $X^*$. Whether $(A_{(1)})_{(1)}=
A_{(1)}?$
\medskip

This question was answered in negative by S.Mazurkiewicz
\cite{Maz}. The result of S.Mazurkiewicz makes it natural to introduce
the following definition. (It was done by S.Banach
\cite{B2, p.~208, 213}. S.Banach used the term ``d\'eriv\'e
faible''.)

\definition{Definition 1.2} For an ordinal $\alpha>1$ the {\it
weak$^*$ sequential closure of order} $\alpha$ of $A$ is the set
$$A_{(\alpha)}=\bigcup_{\beta<\alpha}(A_{(\beta)})_{(1)}.$$
\enddefinition

Weak$^*$ sequential closures were studied by S.Banach in his book \cite{B2}
(see, also, \cite{B3} and \cite{B4}). He proved the following results. 

\proclaim{Theorem 1.3 \cite{B2, p.~124}} Let
$X$ be  a separable Banach space and let $A$ be
a subspace in $X^*$. Then $A=A_{(1)}$ if and only if for
every $f\in X^*\backslash A$ there exists $x\in X$ such that
$f(x)=1$ but $a(x)=0$ for every $a\in A$.
\endproclaim

In modern terminology this result can be stated as: a subspace
in the dual of a separable Banach space is weak$^*$ closed if and only
if it is weak$^*$ sequentially closed.
\medskip 

\proclaim{Theorem 1.4 \cite{B2, p.~213}} Let
$X$ be  a separable Banach space and let $A$ be
a subspace in $X^*$. A necessary and sufficient
condition for $A_{(1)}=X^*$ is that there exists a number $M>0$
such that, for each $x\in X$, the subspace
$A$ contains a functional $f$ satisfying the conditions
$$||f||\le M \ \ \hbox{ and }\ \ |f(x)|=||x||.$$ 
\endproclaim

\proclaim{Theorem 1.5 \cite{B2, pp.~213 and 124}} Let $X$ be
a separable Banach space and let $A$ be a subset in $X^*$.
Then there exists a countable ordinal $\alpha$ such that
$A_{(\alpha)}=A_{(\alpha+1)}$.
\endproclaim

\proclaim{Theorem 1.6 \cite{B2, p.~209}} For every positive integer
$n$ there exists a subspace $A$ in $(c_0)^*=l_1$
such that $A_{(n)}\ne A_{(n+1)}$.
\endproclaim

The book \cite{B2} also contains the following statement:

\proclaim{Statement 1.7 \cite{B2, p.~213}} For every contable ordinal $\alpha$
there exists a subspace $A$ in $(c_0)^*=l_1$
such that $A_{(\alpha)}\ne A_{(\alpha+1)}$.
\endproclaim

This result was not proved in \cite{B2}. At this point S.Banach
referred to a paper that never appeared. Unfortunately the most
comprehensive survey on developements initiated in Banach's
book (I mean the survey by A.Pe\l czy\'nski and Cz.Bessaga
\cite{PB}) does not contain any comments
on this statement. To the best of my knowledge the first proof
of this statement was given by O.C.McGehee \cite{McG}.
\bigskip

It is natural to suppose that the reason for studying
weak$^*$ sequential closures by S.Banach and S.Mazurkiewicz
was the lack of acquaintance of S.Banach and his school
with concepts of general topology.
Although the name ``General topology'' was introduced
later, the subject did already existed, see
F.Hausdorff \cite{H},
Alexandroff-Urysohn \cite{AU} and A.Tychonoff \cite{Ty}.
It is also worth mentioning that J.von Neumann \cite{N, p.~379} had already
introduced the notion of a weak topology.
\medskip

Using the notions of a topological
space and the Tychonoff theorem, more elegant treatment of 
weak and weak$^*$ topologies, and the duality
of Banach spaces was developed
by L.Alaoglu \cite{Al1}, \cite{Al2}, N.Bourbaki \cite{Bou} and 
S.Kakutani \cite{Kak}.
See N.Dunford \& J.T.Schwartz \cite{DS, Sections V.3--V.6} for a well-organized presentation
of this topic.
\medskip

Nevertheless, an ``old-fashioned'' treatment of S.Banach still attracts
attention. It happens because the ``sequential''
approach is very useful in several contexts. Now we would like
to mention some of them.
\medskip

In {\bf Harmonic Analysis} weak$^*$ sequential closures lead to a very useful
and natural classification of sets of uniqueness. This fact was
noticed by I.I.Piatetski-Shapiro (\cite{Pi1}, \cite{Pi2}). 
Since then the classification of sets of uniqueness
was repeatedly used to prove results on sets of
uniqueness using Banach space techniques and results. We are not
going to discuss this topic in any detail.
We refer to A.Kechris and A.Louveau \cite{KL},
A.Kechris, A.Louveau and V.Tardivel \cite{KLT} and
to R.Lyons \cite{Ly}.
\medskip

{\bf Ill-posed problems.} V.A.Vinokurov, Yu.I.Petunin and A.N.Plichko
(\cite{VPP1} and \cite{VPP2}) proved that a very important 
in the theory of ill-posed problems class of regularizable
linear operators can be characterized in terms of weak$^*$ sequential
closures (at least in some important cases).
For further results in this direction see the systematic exposition in
the book of Yu.I.Petunin and A.Plichko \cite{PP} and a more recent paper
\cite{O2}.
\medskip

J.Saint-Raymond \cite{S} and A.N.Plichko \cite{Pl1}, \cite{Pl4}
proved that the {\bf Borel and Baire classification} of inverses of continuous
injective linear operators can be described in terms of weak$^*$
sequential closures (we shall discuss this result in more detail
in \S 3).
\medskip

S.Dierolf and V.B.Moscatelli \cite{DM} found connections
with the {\bf structure theory of Fr\'echet spaces}. See
E.Behrends, S.Dierolf and P.Harmand \cite{BDH}, 
G.Me\-ta\-fune and V.B.Moscatelli \cite{MM1}, \cite{MM2}, 
V.B.Moscatelli \cite{M2} and the present author \cite{O9}
for further results in this direction.
\medskip

A.Plichko \cite{Pl3, Theorem 4} used weak$^*$ sequential
closures to solve a problem on {\bf universal Markushevich
bases in Banach spaces}
posed by N.J.Kalton \cite{Ka, Problem 1, p.~187}.
\bigskip

In connection with these applications of weak$^*$ sequential
closures we find it natural and useful to give an account
on the present state of the study of weak$^*$ sequential closures
initiated by the results of S.Banach and S.Mazurkiewicz.
\medskip
 
Our main purpose is to study the properties of weak$^*$ sequential
closures of total subspaces in the dual spaces of separable Banach
spaces (because it is the case that is the most important for
all listed applications).
\medskip

Recall that a subset $M$ of the dual space $X^*$ is called {\it
total} if for every $0\ne x\in X$ there exists $f\in M$ such that
$f(x)\ne 0$.
\bigskip

Before we turn to the main purpose we would like to make some remarks on
non-separable case.
\bigskip

{\bf 1.} It is clear that for reflexive spaces weak$^*$ sequential closures
coincide with weak sequential closures. It turns out that reflexive
spaces are not the only spaces with this property. 
A.Grothen\-dieck \cite{Gr, p. 168} found nonreflexive spaces $X$ such
that every weak$^*$ convergent sequence in $X^*$ is weakly
convergent. (Now such spaces $X$ are called {\it Grothendieck spaces}.)
It is easy to see that nonreflexive Grothendieck spaces
are non-separable.
See J.Diestel and J.J.Uhl, Jr. \cite{DU, p.~179}
for a survey on Grothendieck spaces and 
J.Bourgain \cite{Bo}, R.Haydon \cite{Ha},
S.S.Khurana \cite{Kh} and M.Talagrand \cite{Ta} for more recent results. 
(Additional references can be found using MathSciNet.)
\medskip

It is easy to verify that for convex subsets in the duals of
Grothendieck spaces the weak$^*$ sequential closure coincide with
the norm closure. B.V.Godun \cite{G1,
Proposition 3} observed that the following
converse to this statement is true: if $A=A_{(1)}$ for every
closed subspace in $X^*$, then $X$ is a Grothendieck space.
\bigskip

{\bf 2.} One of the natural and important questions about weak$^*$ sequential
closures is: when is the second dual $X^{**}$ of a Banach space $X$
equal to the weak$^*$ sequential closure of the canonical image
of $X$ in $X^{**}$? For separable
spaces this question was answered by E.Odell and H.P.Rosenthal
\cite{OR}. They proved the following theorem.

\proclaim{Theorem 1.8} The second dual of a separable Banach space $X$
coincide with the weak$^*$ sequential closure of the canonical
image of $X$ if and only if $X$ does not contain a subspace isomorphic
to $l_1$.
\endproclaim

J.Diestel \cite{D, Chapter XIII} presents a proof of Theorem 1.8
with all necessary preliminaries.
The papers R.J.Fleming \cite{Fl},
R.J.Fleming, R.D.McWilliams and J.R.Retherford \cite{FMR},
A.Grothendieck \cite{Gr}, R.D.McWilliams \cite{McW1}, \cite{McW2} and
H.P.Rosenthal \cite{R1}, \cite{R2}, \cite{R3} contain preceding
and relevant results.
\bigskip

\head \S 2 Existence of total subspaces with long chains of strictly
increasing weak$^*$ sequential closures
\endhead

\definition{Definition 2.1} Let $A$ be a subset in $X^*$.
The least ordinal $\alpha$ for
which $A_{(\alpha)}=A_{(\alpha+1)}$ is called the {\it order} of $A$.
(We use the convention $A_{(0)}:=A$. Hence the order of any weak$^*$
closed subset in $X^*$ is equal to $0$.)
\enddefinition

Theorem 1.5 implies that for separable $X$ the order of any subset 
in $X^*$ is a countable ordinal.
\medskip

This section is devoted to the following {\bf question.}
\medskip

Let $X$ be a separable Banach space. What are the possible orders
of total subspaces $M\subset X^*$?
\bigskip

To answer this question we need

\definition{Definition 2.2} A Banach space $X$ is called
{\it quasi-reflexive} if
$\dim (X^{**}/\pi (X))<\infty $, where $\pi :X\to X^{**}$ is
the canonical embedding and $X^{**}/\pi (X)$ is the quotient space.
\enddefinition

It is worth mentioning that 
classical Banach spaces ($L_p,\ H_p,\ C(K)$) are either reflexive
or non-quasi-reflexi\-ve.
The first example of a non-reflexive quasi-reflexive space was
given by R.C.Ja\-mes \cite{Jam}. (This example is discussed in
\cite{LT} (p.~25).) The term ``quasi-reflexive space''
was introduced by P.Civin and B.Yood \cite{CY}, their paper
contains a systematic study of quasi-reflexive spaces.
\medskip

\proclaim{Theorem 2.3}  {\rm (1)}
If $X$ is a non-quasi-reflexive separable 
Banach  space,
then for every countable ordinal $\alpha$ there exists a total
subspace $\Gamma \subset X^*$ of order $\alpha+1$, that is
$$\Gamma\subset\Gamma_{(1)}\subset\Gamma_{(2)}\subset\dots
\subset\Gamma_{(\alpha)}\subset\Gamma_{(\alpha+1)}=X^*,$$
where all inclusions are proper.
\smallskip

{\rm (2)} If $X$ is a quasi-reflexive separable Banach space,
then a total subspace $\Gamma\subset X^*$ is of
order $1$ if $\Gamma\ne X^*$ and of order
$0$ if $\Gamma=X^*$.
\endproclaim

\remark{Remark 2.4} (a) For a separable Banach space $X$
the whole dual space $X^*$ is the only total subspace of order $0$.
A dense subspace $\Gamma\subset X^*$ satisfying $\Gamma\ne X^*$ is a total
subspace of order $1$. Hence in the statement of Theorem 2.3 $\alpha$
may be equal not only to a usual ordinal number, but also to $0$ or $-1$.
\medskip

(b) By Theorem 1.5 the order of a subspace in the dual of a separable
Banach space is a countable ordinal.
B.V.Godun \cite{G2, Lemma 1}
observed that the order of a subspace in the dual of a separable
Banach space cannot be equal to a
countable limit ordinal. (Recall that an ordinal $\gamma$ is said to
be {\it limit} if it cannot be written in the form 
$\alpha+1$. We consider $0$ and $1$ as non-limit ordinals.)
\medskip

(c) Hence Theorem 2.3 gives a complete answer to
the question above. 
\endremark

\remark{Remark 2.5} Using the fact that a subspace $M$
in $X^*$ can be identified in a canonical way with the
total subspace of the dual of the quotient $X/(M^\top)$ and
the fact that quotients of quasi-reflexive spaces are quasi-reflexive
(see \cite{CY}) it is easy to show that any subspace
in the dual of a quasi-reflexive separable Banach
space has order $1$ (if it is not
weak$^*$ closed) and order $0$  (if it is weak$^*$ closed).
\endremark
\medskip

The first statement of Theorem 2.3 was proved by the present author
in \cite{O1}, the second was known earlier (see below).
Theorem 2.3 (1) has a rather long history and many partial results
were known prior the publication of \cite{O1}.
The history started with
the results of S.Mazurkiewicz and S.Banach mentioned in \S 1
(see Theorems 1.5, 1.6 and Statement 1.7). To the best of my knowledge
the Banach's proof of Statement 1.7 did not exist and
the first example of a Banach space whose dual contain subspaces of
any non-limit countable order was given by
D.Sarason (see \cite{S1} and \cite{S3}). In these papers he
constructed such subspaces in $H_\infty$. Soon afterwards
O.C.McGehee \cite{McG} proved that the same is true for
$X=c_0$ thus proving the Banach's Statement 1.7.
In his proof he identifies $l_1$ with the space of absolutely
convergent Fourier series and uses a quite nontrivial amount
of Harmonic Analysis. A bit later D.Sarason \cite{S2} 
found total subspaces of all possible orders in $l_\infty$.
D.Sarason's proofs in all of the mentioned paper use some tools from Complex
Analysis. The subspaces constructed by O.C.McGehee \cite{McG} and
D.Sarason \cite{S1} are also total. 

\remark{Remark 2.6} Both O.C.McGehee and D.Sarason consider complex spaces.
But their results imply similar results for the real spaces
$c_0$ and $l_1$ also. We discuss this observation in the case $X=c_0$,
the case
$X=l_1$ can be considered in the same way. Observe that the complex
space $c_0$ considered as a real space is isomorphic to the
real space $c_0$.
Let $M$ be the total subspace in the dual of complex $c_0$ of order $\gamma$
constructed by O.C.McGehee. Each $f\in M$ has a unique representation
as $f=f_1+if_2$, where $f_1$ and $f_2$ are real-valued $\Bbb R-$linear
functionals. The set $L:=\{f_1:\ f\in M\}$ is a linear subspace of
the dual of complex $c_0$ considered as a real space.
Using the well-known trick due to Bohnenblust-Sobczyk-Soukhomlinoff
(see \cite{DS, pp.~63--64}) or the explicit
representation of $L$ it is easy to show that $L$ is a total
subspace in the dual of $c_0$ considered as a real space and the order
of $L$ coincide with the order of $M$.
\endremark
\medskip

J.Dixmier \cite{Di} introduced a very useful in the present
context notion of a characteristic
of a subspace in a dual Banach space and found several equivalent
definitions of it. One of the equivalent definitions is:

\definition{Definition 2.7} Let $M$ be a subspace in $X^*$. The
{\it characteristic} of $M$ is defined by
$$r(M):=\inf_{0\ne x\in X}
\frac{\sup\{|f(x)|:\ f\in M,\ ||f||\le 1\}}{||x||}.$$
\enddefinition

It is easy to see that Theorems 1.3 and 1.4 imply the following
result.

\proclaim{Corollary 2.8} The characteristic of a
total subspace in the dual of a separable
Banach space is $>0$ if and only if the order of the subspace
is $\le 1$.
\endproclaim

J.Dixmier observed that the construction of S.Mazurkiewicz gives
an example of a total subspace in $c_0^*$ with zero characteristic.
He constructed a similar example in $l_1^*$
\cite{Di, pp.~1067--1068}.
(His construction works in
more general situation, but the terminology needed to
state more general result was not available at that time.)
\medskip

J.Dixmier's construction is based on the following result.

\proclaim{Theorem 2.9 \cite{Di, p.~1064}} Let $X$ be a Banach space
and let $\Gamma$ be a subspace in $X^*$. Then
$$r(\Gamma)=\inf\left\{\frac{||x+z||}{||x||}:\ x\in\pi(X), x\ne 0,
z\in \Gamma^\bot\right\},$$
where $\Gamma^\bot:=\{z\in X^{**}:\ z(x^*)=0\hbox{ for every }
x^*\in\Gamma\}$.
\endproclaim

Joining this result with Banach's Theorem 1.4 we get
the second statement of Theorem 2.3.
\smallskip

In fact, let $X$ be quasi-reflexive and $\Gamma\subset X^*$ be total.
By definition of a total subspace it implies
$\Gamma^\bot\cap\pi(X)=\{0\}$. The definition of a quasi-reflexive
space immediately implies that $\Gamma^\bot$ is finite-dimensional.
Using Theorem 2.9 and the compactness argument we get
$r(\Gamma)>0$. Jointly with Corollary 2.8 it proves the second part of
Theorem 2.3.
\medskip

The study of total subspaces with zero characteristic was
continued by Yu.I.Petu\-nin \cite{P} who
found a wide class of Banach spaces
whose duals contain total subspaces of characteristic zero,
observed that it cannot happen for quasireflexive spaces and
asked whether such subspaces exist for every non-quasi-reflexive space.
W.J.Davis and J.Lindenstrauss \cite{DL} answered this question in the
affirmative. Joining their result with the mentioned
result of Yu.I.Petunin \cite{P} we get the following
theorem.

\proclaim{Theorem 2.10} 
$X^*$ contains a total subspace with characteristic
zero if and only if $X$ is non-quasi-reflexive.
\endproclaim

A.Plichko \cite{Pl5} found another proof of this result, see
also \cite{Pl2}.
\medskip

W.J.Davis and W.B.Johnson \cite{DJ} found a very important
characterization of non-quasi-reflexive spaces.

\proclaim{Theorem 2.11 \cite{DJ, p.~362}}  A  Banach space $X$ is
non-quasi-reflexive if and only if $X$
contains a bounded away from $0$
basic sequence $\{z_n\}^\infty_{n=0}$ such that the set
$$
\{||\sum^k_{i=j}z_{i(i+1)/2+j}||\}^\infty_{j=0},^\infty_{k=j}
$$
is bounded.
\endproclaim

\remark{Remark 2.12} The difficult part of this theorem is the fact
that each non-quasi-reflexive space contains such basic sequence.
\endremark
\medskip

Using Theorem 2.11 B.V.Godun \cite{G1} proved that for every
$n\in\Bbb N$ the dual of
every non-quasi-reflexive separable Banach space contain a total subspace
$\Gamma$ satisfying $\Gamma_{(n)}\ne\Gamma_{(n+1)}$. His argument can
be used to show that there exists a total subspace $\Gamma$ satisfying
$$\Gamma_{(\omega)}\ne\Gamma_{(\omega+1)},\eqno{(*)}$$
where $\omega$ is the first infinite ordinal.
\medskip

Being unaware of the B.V.Godun's result V.B.Moscatelli (see \cite{Op})
posed a problem
on existence of subspaces satisfying (*) and rediscovered the B.V.Godun's
result (with a bit different proof) in \cite{M1}.
\medskip

B.V.Godun \cite{G2} made an attempt to prove similar results
for larger cardinals, but the inductive argument in \cite{G2} does not work
for infinite ordinals.
\medskip

The final step in proving the first statement of Theorem 2.3 was
made by the present author in \cite{O1}. The argument of \cite{O1} gives
a new proof of the result even for $c_0$. Reading \cite{O1} it may
be useful to look at the special case of $c_0$ first. The approach of
\cite{O1} is such that using the W.J.Davis-W.B.Johnson characterization of non-quasi-reflexive spaces (Theorem 2.11),
the proof for $c_0$ can be easily transferred to any
non-quasi-reflexive space. (It is far from being clear whether the
McGehee's proof can be transferred.)
\smallskip

The paper \cite{O1} is available in English, Russian and
Ukrainian. Since no simplifications of the original
proof have been found since the publication of \cite{O1}
in 1987 I have decided not to reproduce it here.
\medskip

As we have already observed Theorem 2.3 gives a complete answer to 
the question posed at the begining of \S 2. But it is not the end of
the story. The present author generalized Theorem 2.3 in two different
directions. One of the directions is discussed in \S 3 (see
Theorem 3.8). The other can be found in \cite{O7}. 
The present author \cite{O3} studied the isomorphic
structure of total subspaces of given order. More recently,
results of this type attracted attention of A.J.Humphreys and S.G.Simpson
\cite{HS}, whose main interest was foundations of mathematics.
They found another proof of the result of O.C.McGehee. Their
proof does not use Harmonic Analysis. They use the notion
of a well-founded tree. Their proof is a kind of a dual of the
proof in \cite{O1} (if we apply the latter to
$X=c_0$). The use of well-founded trees is quite natural in this
context (and is implicit in \cite{O1}). Nevertheless I do not
feel that their proof is more elementary than the proof in
\cite{O1}.
\medskip

It seems that the question posed at the begining of this section has not been
studied in any detail for general Banach spaces. It is known
that Theorems 1.3--1.5 are not valid for non-separable spaces.
J.Dixmier \cite{Di} found natural and useful analogues of
Theorems 1.3 and 1.4 for general Banach spaces;
these analogues are not in terms of
weak$^*$ sequential closures.
\medskip

See T.Banakh, T.Dobrowolski and A.Plichko \cite{BDP} and 
V.P.Fonf \cite{Fo} for applications of
Theorem 2.3 in contexts that are not mentioned in this survey.
See D.Sarason \cite{S4} for applications of weak$^*$ sequential closures
to certain problems about the invariant subspaces of normal operators
on a Hilbert space.
\medskip

The papers A.A.Albanese \cite{Alb}, V.B.Moscatelli
\cite{M2}  and the present author \cite{O4},
\cite{O5}, \cite{O6}, \cite{O7} and \cite{O10} contain some more
results on classification of total subspaces of characteristic
zero.
\medskip

The papers J.M.Anderson and J.E.Jayne \cite{AJ}, J.E.Jayne \cite{Jay}
and F.Jellett \cite{Jel} contain generalizations of Sarason's results
in another direction. These authors study not the weak$^*$ sequential
closures, but the sets of pointwise limits of a given space of
functions.
\medskip

\head{\S 3 Borel and Baire classification of linear operators}
\endhead

Let $X, Y$ be separable Banach spaces
(in fact, only the fact that $X$ is separable
is important). Let $T: X\to Y$ be an injective continuous
linear operator. Then $T^{-1}:T(X)\to X$ is a well-defined linear
operator. This operator is discontinuous if for some sequence
$\{x_n\}_{n=1}^\infty\subset X$ we have
$\lim_{n\to\infty}\frac{||Tx_n||}{||x_n||}=0$.
On the other hand the well-known Suslin theorem
(one-to-one continuous image of a Borel subset of a complete
separable metric space is a Borel set, see \cite{Kur, p.~487})
implies that $T^{-1}$ is a Borel map.
In this section we consider the {\bf question}: what is the
Borel class of $T^{-1}$?
\medskip

Let us recall necessary definitions
(see \cite{KL, Chapter IV}, \cite{Kur, \S 30, \S 31}).
Let $X$ be a metric space. Let
$\Cal B$ be the smallest collection of subsets of $X$
that
\smallskip

(a) contains all open subsets;
\smallskip

(b) is closed under the operations of complementation, countable union
and intersection.
\smallskip

Sets from $\Cal B$ are called {\it Borel sets}. 
\medskip

We consider the following hierarchy of Borel sets.
\smallskip

For every countable ordinal $\alpha$ we define multiplicative and
additive classes $\alpha$ in the following way.
\smallskip

(a) $\alpha=0$

- Multiplicative class $0$ is the collection of all closed subsets of $X$.

- Additive class $0$ is the collection of all open subsets of $X$.
\smallskip

(b) For $\alpha\ge 1$

- Multiplicative class $\alpha$ is the collection of all countable
intersections $\cap_{n=1}^\infty P_n$, where $P_n$ belongs to some
additive class $\alpha_n$ with $\alpha_n<\alpha$.

- Additive class $\alpha$ is the collection of all countable
unions $\cup_{n=1}^\infty P_n$, where $P_n$ belongs to some
multiplicative class $\alpha_n$ with $\alpha_n<\alpha$.
\medskip

It is easy to see that each Borel set belongs to some
of the additive (or multiplicative) classes.
\medskip

A map $f:X\to Y$ between metric spaces is said to be of
{\it Borel class} $\alpha$ if the set
$f^{-1}(F)$ is a Borel set of multiplicative class $\alpha$
for every closed subset $F\subset Y$.
\medskip

We recall one more natural classification of maps between
metric spaces \cite{Kur, \S 31.IX}.
\medskip

Let $X$
and $Y$ be metric spaces. The  family  of  {\it analytically
representable} maps from $X$ to $Y$
is defined to be the smallest family of maps from $X$  to $Y$
which contains
\smallskip

1) all continuous maps;
\smallskip

2) the limits of convergent sequences of maps belonging to it.
\medskip

This family is representable as a union
$\cup _{\alpha \in \Omega }\Phi_\alpha$, where $\Omega$ is the
set of all countable ordinals and $\Phi_\alpha$ are defined  in
the following way.
\smallskip

1. The class $\Phi_0$ is the set of all continuous maps.
\smallskip

2. The class $\Phi_\alpha\ (\alpha >0)$ consists of all maps
which are limits
of convergent sequences of maps belonging to
$\cup_{\xi<\alpha}\Phi_\xi$.
\smallskip

The class $\Phi_\alpha$ is called {\it Baire class} $\alpha$.
\medskip

It is known (see S.Banach \cite{B1}, K.Kuratowski \cite{Kur, \S 31},
S.Rolewicz \cite{R})
that if $Y$ is a separable Banach  space,
and $\alpha$ is finite then the Borel class $\alpha$
coincides with the Baire class $\alpha$. If $\alpha$ is infinite
then the Borel class $\alpha+1$ coincide with the Baire class $\alpha$.
\medskip

V.A.Vinokurov \cite{V} proved that for arbitrary Banach space
$Y$  the class of all regularizable maps  from $X$  to $Y$
coincide with the Baire class $1$. (Regularizability is a very
important concept in the theory of ill-posed probelms.
We are not going to define it here, see Yu.I.Petunin and A.N.Plichko
\cite{PP} or V.A.Vinokurov \cite{V}.)
\medskip

See R.D.Mauldin \cite{Mau} for an interesting exposition of results
on Borel sets and Baire classes of real-valued functions.
\bigskip

Now we return to the situation described at the begining of this
section.

\proclaim{Theorem 3.1} Let $X, Y$ be Banach spaces. Let $X$ be separable
and let $T:X\to Y$ be an injective continuous linear operator.
Then the map $T^{-1}:T(X)\to X$
is of Borel class $\alpha$ if and only if the subspace
$T^*Y^*$ is of order $\alpha$ (in the sense of Definition 2.1).
\endproclaim

\remark{Remark 3.2} This result was proved by J.Saint-Raymond
\cite{S, Corollaries 42 and 45}.
Somewhat later A.Pli\-chko independently found a different proof
of it (see the summary \cite{Pl1}). Because of some errors in the
Plichko's proof its publication was delayed. Eventually it was published
in \cite{Pl4}. Writing \cite{Pl4} A.Plichko was already aware of
J.Saint-Raymond's paper \cite{S}. 
\endremark
\medskip

To apply results on the existence of total subspaces of given orders
to show the existence of operators with inverses in given
Borel (Baire) classes we need the following result (folk-lore).

\proclaim{Lemma 3.3} Let $X$ be a separable Banach space, let $M$
be a closed total subspace in $X^*$ and let $Y$ be an infinite-dimensional
Banach space. Then there exists an injective continuous linear
operator $T:X\to Y$ such that $T^*Y^*\subset M$ and
$(T^*Y^*)_{(1)}=M_{(1)}$.
\endproclaim

\demo{Proof} Since $X$ is separable, then $B(X^*)=\{x^*\in X^*:\
||x^*||\le 1\}$ is a metrizable separable topological space
in the weak$^*$ topology. Let $\{m_i\}_{i=1}^\infty$ be
a weak$^*$ dense sequence in $B(M)=\{x^*\in M:\
||x^*||\le 1\}$. Let $\{y_i\}_{i=1}^\infty$ be a
basic sequence in $Y$ satisfying $||y_i||\le 1\
(i\in\Bbb N)$ (see \cite{LT, p.~4}).
\medskip

Let $T:X\to Y$ be defined by
$$Tx=\sum_{i=1}^\infty\frac1{2^i}m_i(x)y_i.$$
It is easy to verify that
$$T^*y^*=\sum_{i=1}^\infty\frac1{2^i}y^*(y_i)m_i.$$
From here it is clear that $T^*Y^*\subset M$ and $m_i\subset T^*Y^*$.
The lemma follows. $\square$
\enddemo

Using Theorem 2.3, Theorem 3.1 and Lemma 3.3 we get

\proclaim{Corollary 3.4} Let $X$ be a separable non-quasi-reflexive
Banach space and let $Y$ be an infinite-dimensional Banach space.
Then for every countable ordinal $\alpha\ge 1$ there exists
an injective continuous linear operator $T:X\to Y$ such that the
map $T^{-1}:T(X)\to X$ is of Borel class $\alpha$, but is not of Borel class
$\beta$ for any $\beta<\alpha$.
\endproclaim

If $X$ and $Y$ are ``concrete'' spaces, it is natural to ask: whether
the operator $T$ can be chosen from some ``natural'' class of operators?
We consider one question of this type, many other can be treated
in a similar way.
\medskip

{\bf Question.} Let $F$
and $G$ be Banach function spaces on [0,1] and $T:F\to G$ be  an  injective
linear integral operator with an analytic kernel. What Borel class
can the map $T^{-1}:T(F)\to F$ belong to?
\medskip

Let us give some clarifications.
\smallskip

1. By an analytic kernel we mean a map $K:[0,1]\times [0,1]\to\Bbb C$
with the following property: For some open subsets $\Gamma$ and $\Delta$
of $\Bbb C$ such that
$[0,1]\subset \Gamma $ and $[0,1]\subset \Delta $, there exists
an analytic continuation of $K$  to
$\Gamma \times \Delta $.
\smallskip

2. If the spaces are real then we  consider  map $K$  taking
real values on $[0,1]\times [0,1]$.
\smallskip

3. We assume that $F$ and $G$ are infinite-dimensional and that $F$
is a separable Banach function space on the closed interval $[0,1]$
continuously and injectively embedded into $L_1(0,1)$.
\medskip

\remark{Remark 3.5} It is easy to see that if $F$ and $G$ are
infinite-dimensional function spaces and there exists an injective integral
operator with an analytic kernel $T:F\to G$, then $G$ should contain an
infinite-dimensional subspace consisting of functions having analytic
continuations to some
open subset $\Gamma\subset\Bbb C$ satisfying $\Gamma\supset [0,1]$. 
Diminishing $\Gamma$ if necessary
we may assume that the functions are bounded on $\Gamma$. Using
the well-known
argument (see \cite{LT, p.~4}) it is easy to see that in such a case $G$
contains a basic sequence $\{y_i\}_{i=1}^\infty$ consisting of functions
satisfying 1) $||y_i||_G\le 1$ and 2)
$\sup_{z\in\Gamma}|\tilde y_i(z)|\le 1$, where $\tilde y_i$ is an
analytic continuation of $y_i$. 
\endremark
\medskip

Let $\Delta$ be a bounded open subset of the complex plane  such  that
$\Delta\supset [0,1]$. Let $f:\Delta\to\Bbb C$ be a bounded
analytic function on $\Delta$.
(If we consider real spaces then we assume in addition that $f$
takes real  values  on  the  real  axis.)  Function $f$  generates  a
continuous functional on $L_1(0,1)$ by means of the formula
$$\langle f,x\rangle=\int^1_0f(t)x(t)dt.$$

Since $F$ is continuously embedded into $L_1(0,1)$  then $f$  generates  a
continuous functional on $F$. Let us denote by $U$  the subspace  of $F^*$
consisting of all functionals  of this type. Observe that $U$
is a total subspace of $F^*$. It  follows
from the following facts: 1) $F$ is injectively embedded into $L_1(0,1)$;
2) the set of functionals  generated  by  polynomials  is a total
subspace of $(L_1(0,1))^*$.
\medskip

We need the following variant of Lemma 3.3. 

\proclaim{Proposition 3.6} Let $\alpha\ge 1$ be a countable ordinal.
If $U$ contains a total subspace $L$ of order $\alpha$ (as a subspace
of $F^*$) and $G$ satisfies the condition from Remark 3.5, then there
exists a linear integral operator $T:F\to G$ with an analytic kernel
such that $T^{-1}:T(F)\to F$ is in Borel class $\alpha$, but is not
in Borel class $\beta$ for any $\beta<\alpha$.
\endproclaim

\demo{Proof} We follow the proof of Lemma 3.3. Let $\{y_i\}_{i=1}^\infty
\subset G$ be the sequence from Remark 3.5
and let $\{m_i\}_{i=1}^\infty$ be a weak$^*$ dense sequence
in $\{l\in L:\ ||l||\le 1\}$. Let
$a_i=\sup_{z\in\Delta}|m_i(z)|$ (recall that $U$ and, hence $L$
consists of functionals generated by bounded analytic functions
on $\Delta$.)
\medskip

Let $T:F\to G$ be defined by
$$Tx=\sum_{i=1}^\infty\frac1{2^ia_i}m_i(x)y_i.$$
(We may and shall assume that none of $m_i$'s is zero and hence
none of $a_i$'s is zero.)
\medskip

Observe that $T$ is an integral operator of the required type.
\medskip

The rest of the proof is as in Lemma 3.3 (the only difference
is that we have $T^*G^*\subset\hbox{cl}(L)$ instead of
$T^*G^*\subset L$, but
it does not cause any problems). $\square$
\enddemo

\remark{Remark 3.7} Proposition 3.6 is taken from \cite{O8}. The
idea goes back to L.D.Meni\-khes \cite{Men}.
\endremark
\medskip

To apply Proposition 3.6 we need a condition under which $U$ contains
a subspace of order $\alpha$. One such condition was found by the present
author in \cite{O8}. As is natural in such degree of generality, it is
a condition for existence of subspaces of large orders only.

\proclaim{Theorem 3.8} Let $F$ be a separable
Banach space and $U$ be a total subspace
in $F^*$. Consider elements of $F$ as functionals on $U$. 
We get a subspace in $U^*$. Denote it by $W$.
If {\rm $\hbox{cl}(W)$} is of infinite codimension in $U^*$,
then for every countable ordinal $\alpha$ there exists a total
subspace $L$  of $U$ whose order is $\ge\alpha$.
\endproclaim

The proof of this theorem uses a ramification of the approach of
\cite{O1} to the proof of Theorem 2.3. It is too technical to
be presented here (see \cite{O7} and \cite{O8}).
\medskip

It is natural to ask: how can we check whether the condition of
Theorem 3.8 is satisfied for a given space?
For many classical spaces the answer follows from the following
observation: if $\hbox{cl}(U)$ contains a subspace isomorphic
to $l_1$, then $U^*$ is non-separable and hence 
$\hbox{cl}(W)$ (it is obviously a separable space)
is of infinite codimension in it. For classical non-reflexive Banach
spaces (e.g. $C(0,1),\ L_1(0,1)$) it is not difficult to find
such isomorphic copies of $l_1$.

\remark{Remark 3.9} As a corollary of Theorem 3.8 and Proposition
3.6 we get
a generalization of the results of L.D.Menikhes \cite{Men}
and A.Plichko \cite{Pl2} (see Remark 3 in \cite{O8}). L.D.Menikhes
\cite{Men} proved that there exists an integral operator  from $C(0,1)$
to $L_2(0,1)$ with infinitely differentiable kernel and nonregularizable
($=$not of the Borel class $1$) inverse. A.N.Plichko \cite{Pl2}
constructed operators with these properties
for wide classes of function spaces. 
\endremark

\newpage
\rightheadtext{REFERENCES}


\widestnumber\key{BMMW}
\Refs

\ref\key AJ
\by J.M.Anderson and J.E.Jayne
\paper The sequential stability index of a function space
\jour Mathematika
\vol 20
\yr 1973
\pages 210--213
\endref

\ref\key Al1
\by L.Alaoglu
\paper Weak topologies of normed linear spaces
\jour Annals of Math. (2)
\vol 41
\pages 252--267
\yr 1940
\endref

\ref\key Al2
\bysame
\paper Weak convergence of linear functionals (abstract)
\jour Bull. Amer. Math. Soc.
\vol 44
\yr 1938
\page 196
\endref

\ref\key Alb
\by A.A.Albanese
\paper On total subspaces in duals of spaces  of type $C(K)$ or $L^1$
\jour Proc. Roy. Irish Acad., Sect. A
\vol 93
\yr 1993
\pages 43--47
\endref

\ref\key A
\by P.S.Aleksandrov
\book Introduction to  the  set  theory  and  general
topology
\publ ``Nauka''\publaddr Mos\-cow
\yr 1977
\lang Russian
\transl\nofrills German translation:\ \
\by P.S.Alexandroff
\book Einf\"uhrung in die Mengenlehre und in die allgemeine
Topologie 
\bookinfo Hochschulbücher f\"ur
Mathematik, 85
\publ VEB Deutscher Verlag der Wissenschaften
\publaddr Berlin
\yr 1984
\endref

\ref\key AU
\by P.Alexandroff and P.Urysohn
\paper Zur Theorie der topologischen R\"aume
\jour Math. Annalen
\vol 92
\pages 258--266
\yr 1924
\endref

\ref\key B1
\by S.Banach
\paper \"Uber analytisch darstellbare Operationen in abstakten R\"aumen
\jour Fund.\newline Math.
\vol 17
\yr 1931
\pages 283--295; Reprinted with commentary in: S.Banach,
Oeuvres, v.~I, Warsaw,  PWN-\'Editions Scientifiques de Pologne,
1967, pp.~207--217
\endref

\ref\key B2
\bysame
\book Th\'eorie des op\'erations li\'neaires
\publ Monografje Matematyczne
\publaddr Warszawa
\bookinfo (This edition was reprinted by Chelsea Publishing Company,
without proper reference. It was also reprinted in \cite{B5}) 
\yr 1932
\endref

\ref\key B3
\bysame
\book Kurs funktsional'nogo analizu (A course of functional
analysis)
\lang Ukrain\-ian
\publaddr Kyiv
\publ Radyans'ka Shkola
\yr 1948
\bookinfo (A Ukrain\-ian translation of \cite{B2}, with some alterations.
Alterations are not commented and it is not clear who made them)
\endref

\ref\key B4
\bysame
\book Theory of linear operations
\publ North-Holland
\publaddr Amsterdam New York Oxford Tokyo
\yr 1987
\endref

{\bf Remark on \cite{B4}.} The English translation of the Banach's book
was published when many new results on weak$^*$ sequential closures were
available. Since these results are closely connected to the results
of the appendix to \cite{B4} it was natural to mention them in the
English translation. Unfortunately it was not done.
Also all Banach's footnotes and, in particular, majority of his
references were removed. So the history of the subject cannot be properly
understood without looking at \cite{B2} or \cite{B3}.
\medskip

\ref\key B5
\by S.Banach
\book Oeuvres
\vol II
\publ PWN-\'Editions Scientifiques de Pologne
\publaddr Warsaw
\yr 1979
\endref

\ref\key BDP
\by T.Banakh, T.Dobrowolski and A.Plichko
\paper The topological and Borel classification of operator images
\jour Dissert. Math.
\vol 387
\yr 2000
\pages 37--52
\endref

\ref\key BDH
\by E.Behrends, S.Dierolf and P.Harmand
\paper On a problem of Bellenot and Dubinsky
\jour Math. Ann.
\vol 275
\yr 1986
\pages 337-339
\endref

\ref\key Bou
\by N.Bourbaki
\paper Sur les espaces de Banach
\jour C. R. Acad. Sci. Paris
\vol 206
\yr 1938
\pages 1701--1704
\endref

\ref\key Bo
\by J.Bourgain
\paper $H^\infty$ is a Grothendieck space
\jour Studia Math.
\vol 75 
\yr 1983
\pages 193--216
\endref

\ref\key CY
\by P.Civin and B.Yood
\paper Quasi-reflexive spaces
\jour Proc. Amer. Math. Soc.
\vol 8
\yr 1957
\pages 906-911
\endref

\ref\key DJ
\by W.J.Davis and W.B.Johnson
\paper Basic  sequences  and  norming
subspaces in  non-quasi-reflexive  Banach  spaces
\jour Israel  J. Math.
\vol 14
\yr 1973
\pages 353--367
\endref

\ref\key DL
\by W.J.Davis and J.Lindenstrauss
\paper On total nonnorming subspaces
\jour Proc. Amer. Math. Soc.
\vol 31
\yr 1972
\pages 109--111
\endref

\ref\key DM
\by S.Dierolf and V.B.Moscatelli
\paper A note on quojections
\jour Functiones et approximatio
\vol 17
\yr 1987
\pages 131--138
\endref

\ref\key D
\by J.Diestel
\book Sequences and series in Banach spaces
\bookinfo Graduate Texts in Mathematics, 92
\publ Springer-Verlag
\publaddr New York-Berlin
\yr 1984
\endref

\ref\key DU
\by J.Diestel and J.J.Uhl, Jr.
\book Vector measures
\bookinfo Mathematical Surveys, No. 15
\publ American Mathematical Society
\publaddr Providence, R.I.
\yr 1977
\endref

\ref\key Di
\by J.Dixmier
\paper Sur un th\'eor\`eme de Banach
\jour Duke Math. J.
\vol 15
\yr 1948
\pages 1057--1071
\endref

\ref
\key DS
\by N. Dunford and J.T. Schwartz
\book Linear Operators: General Theory
\bookinfo Pure and Applied Mathematics, vol.7
\publ Interscience
\publaddr New York
\yr 1958
\endref

\ref\key Fl
\by R.J.Fleming
\paper Weak$^*$-sequential closure and characteristic of subspaces
of conjugate Banach spaces
\jour Studia Math.
\yr 1966
\vol 26
\pages 307--313
\endref

\ref\key FMR
\by R.J.Fleming, R.D.McWilliams and J.R.Retherford
\paper On $w^*-$sequential convergence, type $P^*$ bases and
reflexivity
\jour Studia Math.
\yr 1965
\vol 25
\pages 325--332
\endref

\ref\key Fo
\by V.P.Fonf
\paper Two theorems on quasireflexive Banach spaces
\jour Ukrain. Mat. Zh.
\vol 42
\yr 1990
\issue 2
\pages 276--279
\lang Russian
\transl\nofrills English transl. in
\jour Ukrainian Math. J.
\vol 42\yr 1990
\issue 2
\pages 245--247
\endref

\ref\key  G1
\by B.V.Godun
\paper Weak$^*$ derived sets of set of linear functionals
\jour Mat. Zametki
\vol 23
\yr 1978
\pages 607--616
\lang Russian
\transl\nofrills English transl. in
\jour Math. Notes
\vol 23
\yr 1978
\pages 333--338
\endref

\ref\key G2
\bysame
\paper Weak$^*$ derived sets of transfinite order of sets of
linear functionals
\jour Sib. Mat. Zh.
\vol 18
\yr1977
\issue 6
\pages 1289--1295
\lang Russian
\transl\nofrills English transl. in
\jour Siberian Math. J.
\vol 18
\yr 1977
\pages 913--917
\endref

\ref\key Gr
\by A.Grothendieck
\paper Sur les applications lin\'eaires faiblement compactes
d'espaces du type $C(K)$
\jour Canadian J. Math. 
\vol 5
\yr 1953
\pages 129--173
\endref

\ref\key H
\by F.Hausdorff
\book Grundz\"uge der Mengenlehre
\publaddr Leipzig
\publ Veit
\yr 1914
\bookinfo (In 1949 this edition was reprinted by Chel\-sea Publishing Company,
without proper reference)
\endref

\ref\key Ha
\by R.Haydon
\paper A nonreflexive Grothendieck space that does not contain $l_\infty$
\jour Israel J. Math.
\vol 40
\yr 1981
\pages 65--73
\endref

\ref\key HS
\by A.J.Humphreys and S.G.Simpson
\paper Separable Banach space theory needs
strong set existence axioms
\jour Trans. Amer. Math. Soc.
\vol 348
\yr 1996
\issue 10
\pages 4231--4255
\endref

\ref\key Jam
\by R.C.James
\paper Bases and reflexivity of Banach spaces
\jour Ann. of Math. (2)
\vol 52
\yr 1950
\pages 518--527
\endref

\ref\key Jay
\by J.E.Jayne
\paper Spaces of Baire functions. I
\jour Ann. Inst. Fourier (Grenoble)
\vol 24
\yr 1974
\issue 4
\pages 47--76
\endref

\ref\key Jel
\by F.Jellett
\paper The sequential stability index and certain spaces of affine functions
\jour Proc. Amer. Math. Soc.
\vol 83
\yr 1981
\issue 1
\pages 36--38
\endref

\ref\key Kak
\by S.Kakutani
\paper Weak topology, bicompact set and the principle of duality
\jour Proc. Imp. Acad. Tokyo
\vol 16
\yr 1940
\pages 63--67; Reprinted in: S.Kakutani, Selected papers, vol. 1,
Boston, Birkh\"auser, 1986, pp. 208--212 
\endref

\ref\key Ka
\by N.J.Kalton
\paper Universal spaces and universal bases in metric linear
spaces
\jour Studia Math.
\vol 61
\yr 1977
\pages 161--191
\endref

\ref\key KL
\by A.S.Kechris and A.Louveau
\book Descriptive set theory and the structure
of sets of uniqueness
\publ Cambridge University Press
\yr 1987
\endref

\ref\key KLT
\by A.Kechris, A.Louveau and V.Tardivel
\paper The class of synthesizable pseudomeasures
\jour Illinois J. Math.
\vol 35
\yr 1991
\pages 107--146
\endref

\ref\key Kh
\by S.S.Khurana
\paper Grothendieck spaces 
\jour Illinois J. Math.
\vol 22
\yr 1978
\pages 79--80
\endref

\ref\key Kur
\by K.Kuratowski
\book Topology
\vol I
\publ Academic Press and PWN
\publaddr New York and Warsaw
\yr 1966
\endref

\ref\key LT
\by J. Lindenstrauss and L. Tzafriri
\book Classical Banach
spaces I, Sequence spaces
\publ Springer-Ver\-lag
\publaddr Berlin
\yr 1977
\endref

\ref\key Ly
\by R.Lyons
\paper A new type of sets of uniqueness
\jour Duke Math. J.
\vol 57
\yr 1988
\pages 431--458
\endref

\ref\key Mau
\by R.D.Mauldin
\paper Baire functions, Borel sets, and ordinary function systems
\jour Advances in Math.
\vol 12
\yr 1974
\pages 418--450
\endref

\ref\key Maz
\by S.Mazurkiewicz
\paper Sur la  d\'eriv\'ee  faible  d'un  ensemble  de
fonctionnelles lin\'eaires
\jour Studia Math.
\vol 2
\yr 1930
\pages 68--71
\endref

\ref\key McG
\by O.C.McGehee
\paper A proof of a statement of Banach about the weak$^*$ topology
\jour Michigan Math. J.
\vol 15
\yr 1968
\pages 135--140
\endref

\ref\key McW1
\by R.D.McWilliams
\paper Iterated $w^*$-sequential closure of a Banach
space in its second conjugate
\jour Proc. Amer. Math. Soc.
\vol 16
\yr 1965
\pages 1195--1199
\endref

\ref\key McW2
\bysame
\paper On certain Banach spaces which are $w^*$-sequentially dense in
their second duals
\jour Duke Math. J. 
\vol 37
\yr 1970
\pages 121--126
\endref


\ref\key Men
\by L.D.Menikhes
\paper On the  regularizability  of  mappings  inverse  to
integral operators
\jour Doklady AN SSSR
\vol 241
\yr 1978
\issue 2
\pages 282--285
\lang Russian
\transl\nofrills English transl. in
\jour Soviet Math. Dokl
\vol 19
\yr 1978
\pages 838--841 (1979)
\endref

\ref\key MM1
\by G.Metafune and V.B.Moscatelli
\paper Quojections and prequojections
\inbook Advances in the Theory of Fr\'echet spaces
\ed T.Terzio\-\v glu
\publ Kluwer Academic Publishers \publaddr Dordrecht
\yr 1989
\pages 235--254
\endref

\ref\key MM2
\bysame
\paper Prequojections and their duals
\inbook Progress in functional analysis
\eds K.D.Bier\-stedt, J.Bo\-net, J.Horvath and M. Maestre
\bookinfo Pe\~niscola, 1990
\publ North-Holland \publaddr Amsterdam
\yr 1992
\pages 215--232
\endref

\ref\key M1
\by V.B.Moscatelli
\paper On  strongly  non-norming  subspaces
\jour Note Mat.
\vol 7
\yr 1987
\pages 311--314
\endref

\ref\key M2
\bysame
\paper Strongly    nonnorming    subspaces    and prequojections
\jour Studia Math.
\vol 95
\yr 1990
\pages 249-\-254
\endref

\ref\key N
\by J. von Neumann
\paper Zur Algebra der Funktionaloperationen und Theorie der
Normalen Operatoren
\jour Math. Annalen
\vol 102
\yr 1929-1930
\pages 370--427; Reprinted in: J.von Neumann, Collected works, v. II,
Pergamon Press, New York, 1961, pp.~86--143
\endref

\ref\key OR
\by E.Odell and H.P.Rosenthal
\paper A double-dual characterization of separable Banach spa\-ces
containing $l^1$
\jour Israel J. Math.
\vol 20
\yr 1975
\pages 375--384
\endref

\ref\key Op
\by Open Problems
\paper Presented at the Ninth Seminar (Poland - GDR)  on
Operator Ideals  and  Geometry  of  Banach  Spaces
\inbook Forschungsergebnisse Friedrich--Schiller--Universit\"at, Jena,
N/87/28
\nofrills\bookinfo Georgental, April, 1986, Communicated by A.Pietsch
\yr 1987
\endref

\ref\key O1
\by M.I.Ostrovskii
\paper $w^*$-derived  sets  of  transfinite  order  of
subspaces of dual Banach spa\-ces
\jour Dokl.
Akad. Nauk Ukrain. SSR
\vol 1987
\issue 10
\pages 9--12
\lang Russian, Ukrainian
\transl\nofrills An English version of this paper can be
found on the web at 
http://front.math.ucdavis.edu/
\endref

\ref\key O2
\bysame
\paper On  the  problem  of  regularizability  of  the
   superpositions of  inverse  linear  operators
\jour Teor.  Funktsii,
   Funktsion. Anal. i Prilozhen.
\vol 55
\yr 1991
\pages 96--100
\lang Russian
\transl\nofrills English transl. in
\jour J.  Soviet
   Math.
\vol 59\yr 1992
\pages 652--655
\endref

\ref\key O3
\bysame
\paper The structure of  total  subspaces 
of  dual  Banach spaces
\jour Teor. Funktsii, Funktsional. Anal. i  Prilozhen.
\vol 58
\yr 1992
\pages 60--69
\lang Russian
\transl\nofrills English transl. in
\jour J. Math. Sci.
\vol 85
\yr 1997
\pages 2188--2193
\endref

\ref\key O4
\bysame
\paper Total subspaces in dual Banach spaces which
are not  norming  over  any  infinite-di\-men\-sio\-nal  subspace
\jour Studia Math.
\vol 105
\yr 1993
\issue 1
\pages 37--49
\endref

\ref\key O5
\bysame
\paper On the classification of total subspaces of dual
   Banach spaces
\jour C. r. Acad. bulg. Sc.
\vol 45
\yr 1992
\issue 7
\pages 9--10
\endref

\ref\key O6
\bysame
\paper Characterizations of  Banach  spaces  which  are
completions with respect to total nonnorming subspaces
\jour Arch. Math.
\vol 60
\yr 1993
\pages 349--358
\endref

\ref\key O7
\bysame
\paper Total subspaces  with  long  chains  of  nowhere
norming weak$^*$  sequential closures
\jour Note Mat.
\yr 1993
\vol 13
\pages 217--227
\endref

\ref\key O8
\bysame
\paper A  note  on  analytical  representability  of
mappings inverse to integral operators
\jour Matematicheskaya Fizika, Analiz i Geometriya
\vol 1
\yr 1994
\issue 3/4
\pages 513--519
\lang Russian
\transl\nofrills MR98k:47064; An English version of this paper can be
found on the web at
http://front.math.ucdavis.edu/
\endref

\ref\key O9
\bysame
\paper On prequojections and their duals
\jour Revista Mat. Univ. Complutense Madrid
\vol 11
\yr 1998
\pages 59--77
\endref

\ref\key O10
\bysame
\paper Completions with respect to total nonnorming subspaces
\jour Matematicheskaya Fizika, Analiz i Geometriya
\vol 6
\yr 1999
\issue 3/4
\pages 317--322
\endref

\ref\key PB
\by A.Pe\l czy\'nski and Cz. Bessaga
\paper Some aspects of the present theory of Banach spaces
\inbook in: S.Banach, Oeuvres
\vol II
\publ PWN-\'Editions Scientifiques de Pologne
\publaddr Warsaw
\yr 1979
\pages 221--302
\transl\nofrills Reprinted in \cite{B4},
\pages 161--237
\endref


\ref\key P
\by Yu.I.Petunin
\paper Conjugate Banach spaces containing  subspaces  of
zero  characteristic
\jour Dokl.  Akad.  Nauk  SSSR
\vol 154
\yr 1964
\pages 527--529
\lang Russian
\transl\nofrills English transl. in
\jour Soviet Math. Dokl.
\vol 5\yr 1964
\pages 131--133
\endref

\ref\key PP
\by Yu.I.Petunin and A.N.Plichko
\book The theory of  characteristic  of
subspaces and its applications
\publ``Vysh\-cha Shko\-la''\publaddr Kiev
\yr 1980
\lang Russian
\endref

\ref\key Pi1
\by I.I.Piatetski-Shapiro
\paper On the problem of the uniqueness of the expansion of a
function in a trigonometric series
\lang Russian
\jour Moskov. Gos. Univ. U\v c. Zap. Mat.
\vol 155(5) 
\yr 1952
\pages 54--72
\endref

\ref\key Pi2
\bysame
\paper Supplement to the work ``On the problem
of uniqueness of expansion of a function in a trigonometric series''
\lang Russian
\jour Moskov.
Gos. Univ. U\v c. Zap. Mat. 
\vol 165(7)
\yr 1954
\pages 79--97
\endref

\ref\key Pl1
\by A.Plichko
\paper Weak$^*$ sequential closures and $B$-measurability of mappings
inverse to linear continuous operators in $WCG$-spaces (summary)
\jour Sibirsk. Mat. Zh.
\vol 22
\yr 1981
\issue 6
\pages 217
\lang Russian
\endref

\ref\key Pl2
\bysame
\paper Non-norming subspaces and  integral  operators  with
non-regularizable inverses
\jour Si\-birsk. Mat. Zh.
\vol 29
\yr 1988
\issue 4
\pages 208-211
\lang Russian
\transl\nofrills English transl. in
\jour Siberian Math. J.
\vol 29\yr 1988
\pages 687--689
\endref

\ref\key Pl3
\bysame
\paper On bounded biorthogonal systems in  some  function
spaces
\jour Studia Math.
\vol 84
\yr 1986
\pages 25--37
\endref

\ref\key Pl4
\bysame
\paper Decomposition of Banach space into a direct sum of
separable and reflexive subspaces and Borel maps
\jour Serdica Math. J.
\vol 23
\yr 1997
\pages 335-350
\endref

\ref\key Pl5
\bysame
\paper A criterion for the quasireflexivity of a Banach space
\lang Ukrainian
\jour Dopovidi Akad. Nauk Ukrain. RSR Ser. A,
\vol 1974
\pages 406--408
\endref

\ref\key R
\by S.Rolewicz
\paper On  the  inversion  of  non-linear  transformations
\jour Studia Math.
\vol 17
\yr 1958
\pages 79--83
\endref

\ref\key R1
\by H.P.Rosenthal
\paper A characterization of Banach spaces containing $l_1$
\jour Proc. Nat. Acad. Sci. U.S.A.
\vol 71
\yr 1974
\pages 2411--2413
\endref

\ref\key R2
\bysame
\paper Some recent discoveries in the isomorphic theory of Banach spaces
\jour Bull. Amer. Math. Soc.
\vol 84
\yr 1978
\pages 803--831
\endref

\ref\key R3
\bysame
\paper  A characterization of Banach spaces containing $c_0$
\jour J. Amer. Math. Soc.
\vol 7
\yr 1994
\pages 707--748
\endref

\ref\key S
\by J.Saint-Raymond
\paper Espaces a mod\`ele s\'eparable
\jour Annales Inst. Fourier (Grenoble)
\yr 1976
\vol 26
\issue 3
\pages 211--256
\endref

\ref\key S1
\by D.Sarason
\paper On the order of a simply connected domain
\jour Michigan Math. J.
\vol 15
\yr 1968
\pages 129--133
\endref

\ref\key S2
\bysame
\paper A remark on the weak-star topology of $l^\infty$
\jour Studia Math.
\vol 30
\yr 1968
\pages 355--359
\endref

\ref\key S3
\bysame
\paper Weak-star generators of $H^\infty$
\jour Pacific J. Math.
\vol 17
\yr 1966
\pages 519--528
\endref

\ref\key S4
\bysame
\paper Weak-star density of polynomials
\jour J. Reine Angew. Math.
\vol 252
\yr 1972
\pages 1--15
\endref

\ref\key Ta
\by M.Talagrand
\paper Un nouveau ${\Cal C}(K)$ qui poss\`ede la propri\'et\'e
de Grothendieck
\jour Israel J. Math.
\vol 37
\yr 1980
\pages 181--191
\endref

\ref\key Ty
\by A.Tychonoff
\paper \"Uber die topologische Erweiterung von R\"aumen
\jour Math. Annalen
\vol 102
\yr 1929
\pages 544-561
\endref

\ref\key V
\by V.A.Vinokurov
\paper Regularizability and analytical representability
\jour Doklady AN SSSR
\vol 220
\yr 1975
\issue 2
\pages 269--272
\lang Russian
\transl\nofrills English transl. in
\jour Soviet Math. Dokl.
\vol 16\yr 1975
\pages 52--55
\endref

\ref\key VPP1
\by V.A.Vinokurov, Yu.I.Petunin and A.N.Plichko
\paper Conditions for measurability and regularizability
of mappings inverse to linear continuous mappings
\jour Doklady AN SSSR
\yr 1975 \vol 220 \issue 3
\pages 509--511
\lang Russian
\transl\nofrills English transl. in
\jour Soviet Math. Dokl.
\vol 16\yr 1975
\pages 97--99
\endref

\ref\key VPP2
\bysame
\paper Measurability and regularizability of mappings that are
inverses of continuous linear operators
\jour Mat. Zametki
\vol 26
\yr 1979
\issue 4
\pages 583--591
\lang Russian
\transl\nofrills English transl. in
\jour Math. Notes
\vol 26\yr 1979
\pages 781--785
\endref

\endRefs
\enddocument